\documentclass[12pt,reqno,letterpaper]{amsart}
\usepackage{euscript,verbatim}
\usepackage[paper=letterpaper,margin=1in]{geometry}

\newcommand{\bR}{{\mathbb R}}

\title{Gauss Curvature Flow on Surfaces of Revolution - The Noncompact Case} 
\author{Thalia D. Jeffres and Leonardo Solanilla } 
\date{December 13, 2023} 

\begin{document} 
\maketitle

\section{ABSTRACT:} 

Earlier work of the first author [15] examined two boundary value problems associated to the Gauss Curvature 
Flow on a surface of revolution generated by a positive, differentiable function on a compact 
interval. In this continuation, two noncompact cases are addressed. 

{\bf Keywords:} Gauss curvature flow, quasilinear parabolic equations.

\section{Introduction.}

This paper is intended as a continuation of, or postscript to, an earlier paper, [15], in which were 
considered two boundary value problems associated to the Gauss Curvature Flow on a surface of 
revolution generated by a positive differentiable function on a compact interval. Here, we 
consider the noncompact case in which the flow is initiated by the surface of revolution generated 
by a positive differentiable function  \( f: \mathbb{R} \rightarrow \mathbb{R}  \) which is 
bounded from below by a positive constant. We assume 
\( f\in C^{2} (\mathbb{R} ) \) and consider two problems. 

{\bf Problem One.} The initial value function $f$ belongs to \( C^{2} (\mathbb{R} ), \)  with 
both $f$ and its first derivative bounded. 

For this problem, we have the following resolution. 

{\bf Theorem One.} For \( f\in C^{2} (\mathbb{R} ) \) with $f$ and \( f'\) bounded, and $f$ 
bounded below by a positive constant, and for \( T >0 \) and \( \alpha \in (0,1), \) 
the equation 
\[ -\frac{\partial u}{\partial t} + \frac{1}{u (1+u_{x} ^{2} )^{3/2} } \frac{\partial ^{2} u}{\partial x^{2} } =0 \] 
in \( \mathbb{R} \times (0,T) \) has a solution \( u\in H_{2+\alpha } (\mathbb{R} \times (0,T) ) ,\) 
satisfying the initial value \( u(x,0) = f(x) .\)

{\bf Problem Two.} The initial value function $f$ belongs to \( C^{2} (\mathbb{R} ), \) and outside of a compact set, \(  f\equiv c,\) a positive constant. 

{\bf Theorem Two.} Let  \( f \in  C^{2} (\mathbb{R} ),  \) with \( f\equiv c >0 \) outside of a compact set. Then the same equation given in Theorem One has a solution 
\( u \in H_{2+\alpha } (\mathbb{R} \times (0,\infty ) ).\) For each fixed \( t >0,\) solutions \( u(x,t) \rightarrow c \) 
as \( \mid x\mid \rightarrow \infty ,\) and for a time interval \( (0,T) \) for \( T<\infty \) 
and \( \varepsilon >0, \) there exists \( R >0\) such that \( \mid u(x,t) -c \mid < \varepsilon \) 
for all \( (x,t) \) with \( \mid x\mid >R \) and \( t\in (0,T). \) 

A family of embedded surfaces in \( \bR ^{3} \) is said to move under Gauss Curvature Flow if each point 
\( p \) moves in the direction normal to the surface and at a rate proportional to the Gauss curvature at 
that point, that is if each point $p$  satisfies  
\[ \frac{\partial p}{\partial t} = \kappa \cdot n , \]
where $n$ is a chosen unit normal vector field and $\kappa $ is the Gauss curvature function. 
In the case of a closed, convex hypersurface, without boundary, the unit outward normal is 
usually taken, in which case the equation is written as 
\[ \frac{\partial p}{\partial t} = - \kappa \cdot n. \] 

For a point of view that emphasizes classical submanifold theory, one may consult works of 
the second author, [22] and [23]. 

If the surface is a surface of revolution generated at each time $t$ by the curve which 
is the graph of a positive function 
\( u(x,t), \) then the generating function satisfies the parabolic equation of one spatial 
variable given by 
\[ \frac{\partial u}{\partial t} = \frac{1}{u (1+ u_{x} ^{2} ) ^{3/2} } \cdot \frac{\partial ^{2} u}{\partial x^{2} } .\] 
Because $u$ is positive, this equation is parabolic regardless of the sign of \( \partial ^{2} u/\partial x^{2} .\)

Firey [10] initiated the study of the Gauss Curvature Flow, deriving the equation that captures the 
wearing of a smooth stone, idealized to a convex body of uniform density. Assuming existence, 
Firey proved that solutions remain convex, with exponentially decreasing volume, and he conjectured  (and proved, under additional hypotheses) that the shape becomes round. He also established 
directions and tools which proved to be extremely fruitful to later investigators. Early on, 
existence on a finite time interval was established by Tso, [24], for the original equation, and 
by Chow, [8], for hypersurfaces moving under powers of the Gauss curvature, whereupon attention 
turned to understanding the shape of solution hypersurfaces as time approaches the moment at 
which solutions cease to exist by collapsing to a point. Contributions include those of Tso, 
[24], Chow, [8], Andrews [1] and [2], Andrews, Guan, and Ni, [3]; the final piece of the puzzle 
was supplied in 2017 by Brendle, Choi, and Daskalopoulos, [7]. 

Closer to the problem under consideration here is the case of a graph, and here important contributions include those of Oliker, [20], whose results were more recently extended to more 
general powers by Li and Wang, [17], and by Ivochkina and Ladyzhenskaya, [13] and [14]. 

Ishimura, [12], studied self-similar solutions to the same equation of this paper; the self-similarity condition is the same as that which produces the Barenblatt solution to the porous 
medium equation (see [9]). Urbas, [25], also investigated certain self-similar, convex, solutions. In 2014, McCoy, Mofarreh, and Williams, [19], considered axially symmetric boundary value 
problems more general than those of [15]. Barrett, Garcke, and  N\"{u}rnberg, [4], took a 
numerical point of view to axially symmetric curvature flows. 

Along the way, these and other authors observed connections to, and applications to, physical 
processes that go well beyond the wearing of stones that was Firey's original inspiration. These 
include affine geometry, image analysis, grain boundaries, flame fronts, phase boundaries, and the 
movement of biomembranes.

\section{Preliminaries: Notation and Discussion of the Schauder Method.} 

The reader will find it easy to consult the text of Lieberman, [18], from which we use several results, if we adopt 
notations and definitions consistent with that source.  First fixing a positive time $T,$ let \( \Omega  \) denote the 
space-time cylinder,  \( \Omega = \mathbb{R} \times (0, T). \) For a point \( X = (x,t) \in \Omega ,\) 
define 
\[ \mid X \mid = \max \{ \mid x \mid , \, \sqrt{\mid t\mid } \} .\] 
For \( \alpha \in (0,1), \) the parabolic H\"{o}lder space \( H_{2+\alpha } (\Omega ) \) is 
defined to be the set of continuous functions \( u: \Omega \rightarrow \mathbb{R} \) 
possessing derivatives \( (\partial /\partial x)^{k} (\partial /\partial t)^{j} u \) for 
non-negative choices of integers $k$ and $j$ for which \( k+2j\leq 2, \) and having finite 
norm \( \parallel u \parallel _{2+\alpha } \) defined as follows. 
\[ \mid u\mid _{2+\alpha ; \Omega } = \sum _{\beta +2j \leq 2} \sup _{\Omega } \mid D_{x} ^{\beta } D_{t} ^{j} u\mid + \left[ \frac{\partial u}{\partial t} \right] _{\alpha ; \Omega } + \left[ \frac{\partial ^{2} u}{\partial x^{2} } \right] _{\alpha ; \Omega } , \] 
where 
\[ \left[ \frac{\partial u}{\partial t} \right]  _{\alpha ; \Omega }  =  \sup _{X \neq Y} \frac{ \mid \frac{\partial u}{\partial t} (X) - \frac{\partial u}{\partial t} (Y) \mid }{\mid X -Y \mid ^{\alpha } } ,  \] 
and similarly for \( \left[ \partial ^{2} u/\partial x^{2} \right]  _{\alpha ; \Omega } . \) 

We also use the space \( H_{1+\beta } (\Omega ), \) for \( \beta \in (0,1), \) whose norm is 
defined in the same way. That is to say, 
\[ \parallel u\parallel _{1+\beta } = \sum _{k+2j\leq 1} \sup _{\Omega } \mid (\frac{\partial }{\partial x} )^{k} (\frac{\partial u}{\partial t} )^{j} u \mid + \left[ \frac{\partial u}{\partial x} \right] _{\beta ; \Omega } . \] 

The method of solution will be by a fixed point theorem. We record here the version given in 
[9], where it is known as Schaefer's Theorem; other authors (see [11], [18]) regard this statement as a version of, or corollary to, Schauder's Theorem. 

{\bf Theorem (Schaefer's Theorem) [9]}. Let \( \mathcal{B} \)  be a real Banach space, and suppose \( T: \mathcal{B}  \rightarrow \mathcal{B} \) is a continuous and compact mapping. For each \( \sigma \in [0,1], \) associate the operator 
\( \sigma T, \) and consider the union of their fixed points, 
\[ \bigcup  _{\sigma \in [0,1] } \{ u\in \mathcal{B}  \mid (\sigma T) u=u \} . \] 
If this set is bounded in the norm of \( \mathcal{B} ,\)  then the original operator $T$ has a fixed point. 

We apply this theorem with the choice \( \mathcal{B} = H_{1+\beta } (\Omega ). \) This is a 
standard, even ubiquitous, approach, and some basic illustrations can be found in [11] and [18]. 

Until proven 
otherwise, $u$ could vanish, making the equation undefined. To remedy this, let  \( m >0 \) be 
a positive lower bound for $f,$ 
and choose a function \( g(z) \) which is smooth on \( \mathbb{R} ,\)  nondecreasing on \( [0,\infty ) \) and has 
\( g(z) = m/4 \) for  
\( z\leq m/4 ,\) and \( g(z) = z \) for  \( z\geq m/2 .\) This function is then extended by 
symmetry  to nonpositive values of 
$z.$ We apply the existence methods to the modified equation 
\[ -\frac{\partial u}{\partial t} + \frac{1}{g(u) (1+u_{x} ^{2} )^{3/2} } \frac{\partial ^{2} u}{\partial x^{2} } =0, \] 
with the same initial condition \( u(x,0) = f(x). \) It will be shown that solutions fall into the range in 
which \( g(u) = u, \) and therefore any solution to this modified equation in fact satisfies the original equation.

To start the method described in Schaefer's Theorem, take \( v\in H_{1+\beta } (\Omega ) , \) 
for \( \beta \in (0,1), \) and consider the linear initial value problem 
\begin{eqnarray*} 
  -\frac{\partial u}{\partial t} + \frac{1}{g(v) (1+v_{x} ^{2} )^{3/2} } \frac{\partial ^{2} u}{\partial x^{2} } &  =  & 0  \, { \, \rm in \, } \, \Omega , \\ 
  u(x,0) & = & f(x) ,
\end{eqnarray*} 
and let 
\[ T: H_{1+\beta } (\Omega ) \rightarrow H_{1+\beta } (\Omega ) \] 
be the solution operator, \( T(v) = u,\) where $u$ is the solution to this linear initial value 
problem. For 
solvability of the linear equation, see Theorem 9.2.3 of [16]. Solutions $u$ have better regularity, which is why $T$ is compact as a map into \( H_{1+\beta } (\Omega ). \)  Associated to this solution operator $T,$ 
and for each \( \sigma \in [0,1], \) consider the operator \( \sigma T. \) When $u$ 
solves the linear initial value problem above, then \( \sigma u \) solves the same equation in  
\( \Omega ,\) with initial value \( \sigma f \) at \( t=0.\) In this way, existence to the original 
initial value problem is recast as proving a uniform bound on \( \parallel u\parallel _{1+\beta } \) 
for the set of all $u$ which satisfy the initial value problem 
\begin{eqnarray*} 
  -\frac{\partial u}{\partial t} + \frac{1}{g(u) (1+u_{x} ^{2} )^{3/2} } \frac{\partial ^{2} u}{\partial x^{2} } &  =  & 0 \, { \, \rm in \, } \,  \Omega \\ 
  u(x,0) & = &  \sigma \cdot f(x) ,
\end{eqnarray*} 

We are now in a position to apply the method. Each summand in the norm \( \parallel u\parallel _{1+\beta } \) will be bounded separately.

\section{Bounds on \( \mid u\mid .\) } 

In this section, we obtain a bound on \( \sup _{ (x,t) \in  \mathbb{R} \times (0,T)  } \mid u\mid \) for solutions 
to the initial value problem 
\begin{eqnarray*} 
  -\frac{\partial u}{\partial t} + \frac{1}{g(u) (1+u_{x} ^{2} )^{3/2} } \frac{\partial ^{2} u}{\partial x^{2} } &  =  & 0, \,  {\rm in} \, \Omega \\ 
  u(x,0) & = &  \sigma \cdot f(x) ,
\end{eqnarray*} 
for \( \sigma \in [0,1]. \) 

The upper and lower bounds on \( \mid u \mid \) are accomplished by applying a technique explained in [16]. It  requires minor adaptations to this quasilinear 
equation;  because of this, and in order to  make this paper self-contained, all the details are 
presented here. Associated to the modified equation given above, we have the operator $P$ given by 
\[ Pu = - \frac{\partial u}{\partial t} + \frac{1}{g(u) (1+u_{x} ^{2} )^{3/2} } \frac{\partial ^{2} u}{\partial x^{2} } . \] 
 
{\bf Lemma.} Suppose \( u \in H_{1+\beta } (\mathbb{R} \times (0,T) ) \) satisfies the 
differential inequality \( Pu \geq 0 \) and that \( u\leq M \) is known to hold at \( t=0. \) Then 
\( u \leq M \) throughout the region \( \mathbb{R} \times (0,T). \) 

\begin{proof} 
For this fixed choice of $u,$ define another operator \( \overline{P} \) by  
\[ \overline{P} w = -\frac{\partial w}{\partial t} + \frac{1}{g(u) (1+w_{x} ^{2} )^{3/2} } \frac{\partial ^{2} w}{\partial x^{2} } . \] 
Since \( u \in H_{1+ \beta } (\Omega ), \) \( \mid u \mid \) has a finite upper bound \( a \geq 0.\) 
For this $a$ and for each \( R > 0, \) and for \( \lambda > 0 \) to be determined in a moment, 
define a function 
\[ v(x,t) = M + 2 \frac{a}{\cosh R} e^{\lambda t} \cosh x. \] 
Applying the operator \( \overline{P} \) to $v,$ 
\[ \overline{P} v = 2 \frac{a}{\cosh R} e^{\lambda t} \cosh x \left[ -\lambda + \frac{1}{g(u) (1+ (\frac{2a}{\cosh R} e^{\lambda t} \sinh x)^{2} )^{3/2} } \right] . \] 
Recall that \( g(z) \geq m/4 \) for all \( z\in \mathbb{R} , \) where $m$ is a positive lower 
bound for $f.$ 
By choosing 
\[ \lambda > \frac{4}{m} , \] 
then \( \overline{P} v <0 \) in \( \mathbb{R} \times (0,T). \)  Evaluated at $u,$ \( \overline{P} u = Pu \geq 0,\) so \( \overline{P} u - \overline{P} v > 0\) in all of \( \mathbb{R} \times (0,T). \) 

We observe that when \( t=0,\) 
\[ (u-v) (x,0) < 0, \] 
because \( u\leq M \) is known to hold for \( t=0. \) If \( \mid x\mid  = R \) with \( t\in (0,T), \)   then 
\begin{eqnarray*} 
 (u-v) (\pm R ,t) & = & (u(\pm R,t) -ae^{\lambda t} ) -M -ae^{\lambda t} \\ 
      & = & (u(\pm R, t) - ae^{\lambda t} ) - (M +ae^{\lambda t} )  \\ 
      & \leq & 0, 
\end{eqnarray*} 
so long as \( \lambda \geq 0, \) and using the observation that \( M \geq -a. \)

Now let \( \varepsilon \) be small enough that \( T -\varepsilon >0, \) and examine the values of \( u-v \) in the region \( B_{R} (0) \times (0, T-\varepsilon ). \) Might \( u-v \) achieve positive values in the region \( B_{R} (0) \times (0,T-\varepsilon ) \)? 
\( u-v \) is continuous on all of \( \overline{B_{R} (0)} \times [0,T-\varepsilon ], \) a compact 
set, so \( u-v \) attains a maximum in this set. If \( u-v \) has positive values somewhere, 
then this maximum value is positive. Since \( u-v \leq 0 \) when \( t=0 \) or when \( (x,t) \in \{ -R, R \} \times [0,T), \) this maximum occurs at a point \(  X_{0} = (x_{0} ,t_{0} ) \) with \( \mid x_{0} \mid < R \) and \( t\in (0,T-\varepsilon ). \) At such a point, \( t_{0} \leq T-\varepsilon < T, \) 
so \( \partial u/\partial t\) is still continuous here, and therefore 
\[ \frac{\partial }{\partial t} (u-v) (x_{0} ,t_{0} )  \geq 0 \] 
at this point. We also have 
\[ \frac{\partial }{\partial x} (u-v) (X_{0}  ) = 0,  \, {\rm and} \,  \frac{\partial ^{2} }{\partial x^{2} } (u-v)(X_{0} )  \leq 0. \] 
Using the fact that \( u_{x} (X_{0} ) = v_{x} (X_{0} ) \) at the maximum point, we  have 
\[  (\overline{P} u -\overline{P} v) (X_{0} ) = - \frac{\partial }{\partial t} (u-v) (X_{0} ) + \frac{1}{g(u) (1+u_{x} ^{2} )^{3/2} } \frac{\partial ^{2} u}{\partial x^{2} } (u-v) (X_{0} ) \leq 0. \] 

Combined with the earlier calculation, this yields the impossibility 
\[ 0 < (\overline{P} u-\overline{P} v ) (X_{0} ) \leq 0, \] 
from which we conclude that \( u-v \leq 0 \) in \( \overline{B_{R} (0) } \times [0,T-\varepsilon ]. \) Since \( \varepsilon \) is arbitrarily small, this implies that \( u-v \leq 0 \) in \( \overline{B_{R} (0) } \times [0,T). \) Choose any \( (x,t) \in \mathbb{R} \times (0,T). \) By the 
above calculations, 
\[ u(x,t) -M -\frac{2a}{\cosh R} e^{\lambda t} \cosh x \leq 0 \] 
holds for all \( R > \mid x\mid . \) Therefore, \( u(x,t) -M \leq 0. \) This concludes the 
proof of the lemma. 
\end{proof}

In the same way, one may establish that the infemum of the initial value function $f$ serves 
as a lower bound for solutions at all times.

To apply these results to the problem, recall from the previous section that what is needed 
is a uniform bound on all solutions to the family of initial value problems 
\begin{eqnarray*} 
  -\frac{\partial u}{\partial t} + \frac{1}{g(u) (1+u_{x} ^{2} )^{3/2} } \frac{\partial ^{2} u}{\partial x^{2} } &  =  & 0, { \, \rm in \, } \Omega \\ 
  u(x,0) & = &  \sigma \cdot f(x) ,
\end{eqnarray*} 
for \( \sigma \in [0,1]. \) The original initial value function $f$ satisfies 
\[ 0 < m \leq f(x) \leq \sup _{x\in \mathbb{R} } f(x). \] 
This implies that for all \( \sigma \in [0,1], \) 
\[ 0 \leq \sigma f \leq \sup _{x\in \mathbb{R} } f(x), \] 
and therefore $u$ satisfies \( 0 \leq u \leq \sup f. \) 

Lastly, if $u$ solves the equation for \( \sigma =1,\) then $u$ satisfies 
\[ m \leq u \leq \sup f, \] 
and in this range, \( g(u) =u. \) That means that once the steps for the Schauder method are 
carried out, we will obtain a solution $u$ that satisfies the original equation, 
\[ \frac{\partial u}{\partial t} = \frac{1}{u (1+ u_{x} ^{2} ) ^{3/2} } \cdot \frac{\partial ^{2} u}{\partial x^{2} } .\] 

\section{Gradient Bound for Problem One and conclusion of the proof of Theorem One.}

We apply the Bernstein technique. In addition to the original sources, [5] and [6], the paper 
of Serrin, [21] gives a clear expositon of this method in both elliptic and parabolic 
settings.  The starting point for this method is 
the observation that if $u$ solves the heat equation, then 
the combinations \( \parallel \nabla u \parallel, \, 1+ \parallel \nabla u \parallel ^{2} , \) 
and \( \sqrt{1+ \parallel \nabla u \parallel ^{2} } \) are all subsolutions to that same 
equation. Here, \( \nabla u \) means the full gradient, not just the gradient 
in the spatical directions. 

We have already seen that if $u$ satisfies the equation 
\[ - \frac{\partial u}{\partial t} + \frac{1}{g(u) (1+u_{x} ^{2} )^{3/2} } \frac{\partial ^{2} u}{\partial x^{2} } = 0. \] 
then 
\[ m \leq u \leq \sup f, \] 
and and that for these values, \( g(u) =u, \) so this time, we will define the 
auxilliary operator by  \( \overline{P} \) by 
\[ \overline{P} = -\frac{\partial }{\partial t} + \frac{1}{u (1+u_{x} ^{2} )^{3/2} } \frac{\partial ^{2} }{\partial x^{2} } . \] 
Applying  \( \overline{P} \) to the quantity \( v = 1+u_{x} ^{2} \) gives 
\[ \overline{P} v = -\frac{\partial v}{\partial t} + \frac{1}{u (1+u_{x} ^{2} )^{3/2} } \frac{\partial ^{2} v}{\partial x^{2} }. \] 
From the definition of $v,$ 
\[ \frac{\partial v}{\partial t} = 2 \frac{\partial u}{\partial x} \frac{\partial ^{2} u}{\partial t \partial x} , \] 
and from the equation satisfied by $u,$ 
\[ \frac{\partial ^{2} u}{\partial t \partial x} = \frac{1}{u^{2} \left( 1+ u_{x} ^{2} \right) ^{5/2} } \left[ -u_{x} u_{xx}  \left( 1+u_{x} ^{2} \right)  - 3u u_{x} u_{xx} ^{2} + u \left( 1+u_{x} ^{2} \right)  u_{xxx} \right]  . \] 
Then 
\[ \overline{P} v = \frac{2}{u^{2} \left( 1+u_{x}  ^{2} \right) ^{5/2} } \left( u (1+4u_{x} ^{2} ) u_{xx} ^{2} + u_{x} ^{2} \left( 1+u_{x} ^{2} \right)  u_{xx} \right) . \] 
Since \( 1 +4u_{x} ^{2} \geq 1 + u_{x} ^{2} , \) it follows that 
\[ \overline{P} v \geq \frac{2}{u (1+u_{x} ^{2} )^{3/2} } (u_{xx} ^{2} + \frac{u_{x} ^{2} }{u} u_{xx} ). \] 
Completing the square, the factor in parentheses can be rewritten as  
\[ u_{xx} ^{2} + \frac{u_{x} ^{2} }{u} u_{xx} = \left( u_{xx} + \frac{u_{x} ^{2} }{2u} \right) ^{2} - \left( \frac{u_{x} ^{2} }{2u} \right) ^{2} ; \] 
we now have  
\[ \overline{P} v \geq \frac{2}{u(1+u_{x} ^{2} )^{3/2} } \left( -\left( \frac{u_{x} ^{2} }{2u} \right) ^{2}  \right)  = \frac{-u_{x} ^{4} }{2u^{3} \left( 1+u_{x} ^{2} \right) ^{3/2} } . \] 
This implies that 
\[ \overline{P} v >  \frac{-1}{2u^{3} } (1+u_{x} ^{2} )^{1/2} , \] 
or 
\[ \overline{P} v + \frac{1}{2u^{3} } v^{1/2} > 0. \] 
Since \( v = 1+u_{x} ^{2} \geq 1, \, v\geq v^{1/2} , \) and therefore 
\[ \overline{P} v + \frac{1}{2u^{3} } v >0. \] 
If \( w(x,t) = e^{-1/(2M^{3} ) } v(x,t), \) then $w$ satisfies \( \overline{P} w >0. \) Since 
\( \overline{P} \) is a linear operator, one can apply a standard maximum principle (or the 
technique and method from Section Four), to obtain 
\[ w(x,t) \leq \max _{x\in \mathbb{R} } w(x,0) , \] 
which implies that 
\[ (1+u_{x} ^{2} ) (x,t) \leq e^{T/(2M^{3} )} \max (1+ f'^{2} ) . \] 

We can now complete the proof of Theorem One.

{\bf Proof of Theorem One.} 
\begin{proof} With bounds on \( \mid u\mid \) and \( \mid u_{x} \mid \) in place, the remaining 
step is the improvement to the H\"{o}lder bound on the gradient. This is accomplished by 
application of Theorem 12.10 of the text of Lieberman, [18]; see also [14]. 
\end{proof} 

\section{Gradient Bound for Problem Two and Conclusion of the Proof of Theorem Two.}

The first step is to show that \( u \rightarrow c \) as \( \mid x\mid \rightarrow \infty .\) This 
is true for each $t,$ and the convergence is uniform on time intervals of the form \( (0,T). \) 

The proof is accomplished by a comparison argument, trapping the solution beneath a supersolution 
with exponential shape. Before formulating this result precisely, we provide some 
preliminary observations and calculations. 

The prototype is \(  \mathcal{P} (x) = e^{-x^{2} } .\) The maximum value of this function's second 
derivative is 
\[ \mathcal{P} ''(\pm \sqrt{3/2 } ) = 4e^{-3/2} .\] 
We now want to construct a modification that encloses the nonconstant portion of the initial 
value function $f.$ Inserting a parameter $h$ that affects the height, and another parameter $w$ 
that changes the width, define an exponential function 
\[ \mathcal{E} (x) = he^{-wx^{2} } , \] 
for \( h,w >0. \) This function has a maximum height of $h,$ with maximum values of its second 
derivative occurring at \( x = \pm \sqrt{3/2w } .\) At this point, 
\[ \mathcal{E} ''\left( \pm \sqrt{\frac{3}{2w} } \right)  = 4hwe^{-3/2} .\] 

If \( w \in (0,1), \) then \( he^{wx^{2} } \geq he^{ -x^{2}} , \) and this will be used later.

Let \( M = \max f, \) and choose \( R_{0} >0 \) large enough that \( f(x) =c \) for all \( \mid x \mid > R_{0} . \) The portion of the graph of $f$ which is not constant is contained in the rectangle with corners \( (\pm R_{0} ,0) \) and \( (\pm R_{0} ,M). \) The corner point \( (R_{0} ,M) \) will lie below the corresponding point \( (R_{0} ,\mathcal{E} (R_{0} ) ) \) on the graph of \( \mathcal{E} \) if $h$ is chosen large enough that \( he^{-R_{0} ^{2} } > M, \) or \( h> Me^{R_{0} ^{2} } . \) With $w$ chosen in the interval \( (0,1), \) the graph of \( \mathcal{E} \) encloses the rectangle and the graph of $f$ lies below that of \( c+ \mathcal{E} \) because 
\[ c+he^{-w x^{2} } \geq c+he^{-x^{2} } \geq f(x) . \] 
We are now in a position to prove the lemma. 

{\bf Lemma.} Suppose that \( u\in H_{2+\alpha } (\mathbb{R} \times (0,T) ) \) is a solution 
to 
\[ -\frac{\partial u}{\partial t} + \frac{1}{g(u) (1+u_{x} ^{2} )^{3/2} } \frac{\partial ^{2} u}{\partial x^{2} } = 0 \] 
on \( \mathbb{R} \times (0,T), \) satisfying \( u(x,0) = f(x), \) where $f$ is positive everywhere 
and \( f(x) = c \) outside of a compact set. Then for any \( \varepsilon > 0, \) there exists  
\( R > 0 \) such that \( \mid u-c \mid < 3 \varepsilon \) for \( \mid x \mid > R \) and \( t\in (0,T). \)

{\bf Proof of the lemma:} 
\begin{proof} 
Let $u$ be a positive solution to the equation 
\[ -\frac{\partial u}{\partial t} + \frac{1}{g(u) (1+u_{x} ^{2} )^{3/2} } \frac{\partial ^{2} u}{\partial x^{2} } = 0, \] 
satisfying \( u(x,0) = f(x), \) and define the quasilinear operator \( \overline{P} \) by 
\[ \overline{P} w = -\frac{\partial w}{\partial t} + \frac{1}{g(u) (1+w_{x} ^{2} )^{3/2} } \frac{\partial ^{2} w}{\partial x^{2} } .\] 
For arbitrary \( \varepsilon ,\) let 
\[ v(x,t) = c + \varepsilon + \varepsilon \frac{t}{t+1} + he^{-wx^{2} } ; \] 
for now, \( h> Me^{R_{0} ^{2} } \) and \( w\in (0,1), \) but further refinements will be made. Is it 
possible to choose these parameters in such a way that $v$ becomes a strict supersolution? This will be 
the case if 
\[ -\varepsilon \frac{1}{(t+1)^{2} } + \frac{1}{g(u) (1+v_{x} ^{2} )^{3/2} } \frac{\partial ^{2} v}{\partial x^{2} } < 0. \] 
It will be enough if  
\[ \frac{1}{g(u) (1+ v_{x} ^{2} )^{3/2} } \frac{\partial ^{2} v}{\partial x^{2} } < \varepsilon \frac{1}{(T+1)^{2} } .\] 
Whenever  \( \frac{d^{2} }{dx^{2} }  (he^{-wx^{2} } ) \leq 0, \) we have 
\[ \frac{1}{g(u) (1+ v_{x} ^{2} )^{3/2} } \frac{\partial ^{2} v}{\partial x^{2} } \leq 0 < \varepsilon \frac{1}{(T+1)^{2} } , \] 
and when \( \frac{d^{2} }{dx^{2} }  (he^{-wx^{2} } ) \geq 0, \) 
\[ \frac{1}{g(u) (1+v_{x} ^{2} )^{3/2} } \frac{\partial ^{2} v}{\partial x^{2} } \leq \frac{1}{g(u) } \frac{d^{2} }{dx^{2} } (he^{-wx^{2} } ) \leq \frac{4}{m} \frac{d^{2} }{dx^{2} } (he^{-wx^{2} } ) , \] 
remembering that \( g(u) \geq m/4 = (1/4) \min f \) . So, it will be enough to make 
\[ \frac{4}{m} \frac{d^{2} }{dx^{2} } (he^{-wx^{2} } ) < \varepsilon \frac{1}{(T+1)^{2} } .\] 
We already calculated above that 
\[ \frac{d^{2} }{dx^{2} } (he^{-wx^{2} } ) \leq 4w he^{-3/2} , \] 
which means that we want 
\[ 4whe^{-3/2} \leq \frac{m}{4} \varepsilon \frac{1}{(T+1)^{2} } .\] 
For concreteness, choose \( h = 2M e^{R_{0} ^{2} } . \) Then the above inequality holds if \( w \in (0,1) \) and also satisfies 
\[ w < \varepsilon \frac{m}{4} \frac{1}{(T+1)^{2} } \frac{1}{8M} e^{3/2 - R_{0} ^{2} } . \] 
With these choices, \( \overline{P} v <0, \) so that \( \overline{P} u - \overline{P} v = Pu -\overline{P} v >0, \) and this holds at all points in \( \mathbb{R} \times (0,T). \) Also, at  \( t=0, \) 
\[ u(x,0) - v(x,0) = f(x) -(c+2\varepsilon + he^{-wx^{2} } ) < -\varepsilon <0. \] 

We now apply the technique explained in Krylov's book [16], which was also employed earlier, to obtain the bound on \( \mid u\mid .\) 

Since \( u \in H_{2+\alpha } (\mathbb{R} \times (0,T) ), \) there exists a positive number 
$a$ such that \( \mid u(x,t) \mid < a \) for all \( (x,t) \in \mathbb{R} \times (0,T). \) For this 
$a$ and for each \( R >0, \) define 
\[ v_{R} (x,t) = c + \varepsilon + \varepsilon \frac{t}{t+1} + he^{-wx^{2} } + 2\frac{a}{\cosh R} e^{\lambda t} \cosh x. \] 
At \( t=0,\) we have \( u(x,0) < c+ he^{-wx^{2} } < v_{R} (x,0). \) What is true at \( \mid x\mid = R\)? At these points, we have 
\begin{eqnarray*} 
 u(\pm R,t) - v_{R} (\pm R,t) & = & u(\pm R,t) - \left( c+ \varepsilon + \varepsilon \frac{t}{t+1} + he^{-wR^{2} } \right)  - 2ae^{\lambda t} \\ 
  & = & (u(\pm R,t) -ae^{\lambda t} ) - \left( c +\varepsilon + \varepsilon \frac{t}{t+1} + he^{-wR^{2} } + a e^{\lambda t} \right) . 
\end{eqnarray*} 
Since \( \mid u(x,t) \mid \leq a \) for all \( (x,t) \in \mathbb{R} \times (0,T), \) then in the 
first term, 
\[ u(\pm R,t) -a e^{\lambda t} \leq a -a e^{\lambda t} = a(1-e^{\lambda t} ) \leq 0, \] 
so long as \( \lambda \geq 0. \) Therefore,  \( (u-v_{R} )(\pm R,t) \leq 0.\)

We calculate 
that 
\begin{eqnarray*} 
 \overline{P} v_{R} & = & -\varepsilon \frac{1}{(t+1)^{2} } - 2 \lambda \frac{a e^{\lambda t} }{\cosh R} \cosh x \\ 
  & + & \frac{1}{g(u) (1+(v_{R} )_{x} )^{3/2} } \left( \frac{d^{2} }{dx^{2} } he^{-wx^{2} } + 2\frac{a}{\cosh R} e^{\lambda t} \cosh x \right) .  
\end{eqnarray*} 

Grouping the terms, 
\[ \frac{1}{g(u) (1+(v_{R} )_{x} ^{2} )^{3/2} } \frac{d^{2} }{dx^{2} } he^{-wx^{2} } \leq 0 < \varepsilon \frac{1}{(t+1)^{2} } \] 
where \( d^{2} /dx^{2} he^{-wx^{2} } \leq 0, \) and 
\begin{eqnarray*} 
 \frac{1}{g(u) (1+ (v_{R} )_{x} ^{2} )^{3/2} } \frac{d^{2} }{dx^{2} } he^{-wx^{2} } & \leq & \frac{1}{g(u) } \frac{d^{2} }{dx^{2} } he^{-wx^{2} } \\ 
  & \leq & \frac{4}{m} \frac{d^{2} }{dx^{2} } he^{-wx^{2} } \\ 
  & \leq & \varepsilon \frac{1}{(t+1)^{2} } 
\end{eqnarray*} 
when \( d^{2} /dx^{2} he^{-wx^{2} } \geq 0, \) by the choices of $h$ and $w$ already made. In 
the other two terms, 
\[ \frac{1}{g(u) (1+(v_{R} )_{x} ^{2} )^{3/2} } \leq \frac{4}{m} , \] 
and so choosing \( \lambda > 4/m \) makes 
\[ 2 \frac{a e^{\lambda t} }{\cosh R}  \cosh x \left( -\lambda + \frac{1}{g(u) (1+(v_{R} )_{x} ^{2} )^{3/2} } \right)  < 0. \] 
Adding these up, we conclude that \( \overline{P} v_{R} < 0, \) and so \( \overline{P} u - \overline{P} v_{R}  = Pu - \overline{P} v_{R}  >0, \) and this is true in all of \( \mathbb{R} \times (0,T) .\) 

Using the same sign analysis as in the bound on \( \mid u \mid ,\) with the function 
\( v(x,t) \) playing the role of $M,$  it is possible to show that in 
the region \( B_{R} (0) \times (0,T)  ,\) the function \( u-v_{R} \) cannot attain a maximum 
away from the parabolic boundary. This implies that for all \( (x,t) \in B_{R} (0) \times (0,T), \) that 
\[ u(x,t) \leq c+ \varepsilon + \varepsilon \frac{t}{t+1} + he^{-wx^{2} } < c+2\varepsilon + h e^{-wx^{2} } + 2\frac{a}{\cosh R} e^{\lambda t} \cosh x. \] 
Choose \( (x,t) \in \mathbb{R} \times (0,T). \) For all \( R > \mid x\mid ,\) the above inequality 
holds. Therefore, 
\[ u(x,t) \leq c +\varepsilon + \varepsilon \frac{t}{t+1} + he^{-wx^{2} } < c+2\varepsilon + he^{-wx^{2} } . \] 
Since for all \( \varepsilon > 0, \) it is possible to find \( K >0 \) with \( he^{-wx^{2} } < \varepsilon \) for \( \mid x\mid > K \) and \( t\in (0,T), \) this also shows that there exists 
\( K>0\) such that for all \( \mid x\mid  > K \) and \( t\in (0,T), \) that 
\[ u(x,t) \leq c + 3\varepsilon .\] 

By the same reasoning, we can show that for \( \mid x\mid \) large enough, \( u(x,t) >c-3\varepsilon .\) Therefore, \( u \rightarrow c \) as \( \mid x\mid \rightarrow \infty . \) Since the choice of $w$ that works for $T$ also works for all \( t\in (0,T), \) there is a single choice of \( K >0 \) such that 
\[ \mid u(x,t) -c\mid < 3\varepsilon \] 
for all \( (x,t) \in ((-\infty ,-K) \cup (K, \infty ) ) \times (0,T). \) This completes the proof 
of the lemma. 
\end{proof}

{\bf Proof of Theorem Two.} 
\begin{proof} 
The \( C^{0} \) bound is already in place from Section Four. For the 
gradient bound, note that if a function in \( H_{2+\alpha } (\mathbb{R} \times (0,T) ) \) approaches a constant as \( \mid x\mid \rightarrow \infty ,\) then also \( \mid u_{x} \mid \rightarrow 0. \) Choose \( K >0 \) large enough that \( \mid u_{x} \mid < \max \mid f' \mid \) for \( \mid x\mid > K. \) Now apply a standard maximum principle for the gradient on any bounded set \( [-K',K'] \times (0,T). \) (See Theorem 9.7 of [18], and also [15].) In this way, we obtain the improved gradient 
bound 
\[ \min _{x\in \mathbb{R} } f'(x) \leq u_{x} (x,t) \leq \max _{x\in \mathbb{R} } f'(x), \] 
for all \( (x,t) \in \mathbb{R} \times (0,T). \) For the improvement to a H\"{o}lder bound on 
\( u_{x} ,\) we apply Theorem 12.10 of [18]. By the Schauder method, the solution exists on \( \mathbb{R} \times (0,T), \) and in fact on \( \mathbb{R} \times (0, \infty ) , \) since $T$ is 
arbitrary and the estimates do not involve $T.$ 
\end{proof}

\vspace{2cm}

\noindent {\bf Acknowledgments} It was while working at the Universidad Michoacana that we 
discovered a mutual interest in Gauss curvature. We thank our former colleagues at the Instituto de F\'{i}sica y Matem\'{a}ticas of that university. We also gratefully acknowledge a conversation 
with Lu Peng. 

\section{References} 

[1] B. Andrews, {\em Motion of hypersurfaces by Gauss Curvature,} Pacific J. Math. 195, 1-34 (2000)

[2] B. Andrews, {\em Gauss Curvature flow: the fate of the rolling stones,} Invent. Math. 138, 151-161 (1999)

[3] B. Andrews, P. Guan, and L. Ni, {\em Flow by powers of the Gauss curvature,} Adv. Math. 299 (2016), 174-201. 

[4] J. Barrett, H.  Garcke, and R.  N\"{u}rnberg,   {\em Variational discretization of axisymmetric curvature flows,}  Numerische Mathematik (2019) 141:791-837. 

[5] S. Bernstein, {\em Sur la g\'{e}n\'{e}ralisation du probl\`{e}me de Dirichlet II,} Math. Ann. {\bf 69} (1910), 82-136. 

[6] S. Bernstein, {\em Sur les \'{e}quations du calcul des variations,} Ann. Sci. \'{E}cole Norm. 
Sup. {\bf 29} (1912), 431-485. 

[7] S. Brendle, K. Choi, and P. Daskalopoulos, {\em Asymptotic Behavior of Flows by Powers of the Gaussian Curvature,} Acta Math. 219(1): 1-16 (September 2017).

[8] B. Chow, {\em Deforming convex hypersurfaces by the $n$th root of the Gaussian curvature,} J. Differential Geom. 22 (1995), no. 1, 117-138. 

[9] L.C. Evans, {\em Partial Differential Equations,} American Mathematical Society, Ed. 2, Providence, 2010. 

[10] W. Firey, {\em Shapes of Worn Stones,} Mathematika, A Journal of Pure and Applied 
Mathematics, Vol. 21, Part 1, No. 41, June, 1974. 

[11] D. Gilbarg and N.S. Trudinger, {\em Elliptc Differential Equations of Second Order,} Springer-Verlag, Berlin, 2001. Reprint of the third edition (1998). 

[12] N. Ishimura, {\em Self-Similar Solutions for the Gauss Curvature Evolution of Rotationally 
Symmetric Surfaces,}  Nonlinear Anal Theory Methods and Applications {\bf 33} no. 1, (1998), 97-104. 

[13] N.M. Ivochkina and O.A. Ladyzhenskaya, {\em Flows generated by symmetric functions of the 
eigenvalues of the Hessian,} Journal of Mathematical Sciences, Vol. 87, No. 2, 1997.

[14] N.M. Ivochkina and O.A. Ladyzhenskaya, {\em On parabolic equations generated by symmetric functions of 
the principal curvatures of the evolving surface, or of the eigenvalues of the Hessian, Part I: 
Monge-Amp\`{e}re Equations,} St. Petersburg Math. J. Vol. 6 (1995), No.3. 

[15] T. Jeffres, {\em Gauss Curvature Flow on Surfaces of Revolution,}  Advances in Geometry, {\bf 9} (2009), no. 2, 189-197. 

[16] N.V. Krylov, {\em Lectures on Elliptic and Parabolic Equations in H\"{o}lder Spaces,} 
American Mathematical Society, Providence, Rhode Island, 1996. Reprint of the first edition (1997).

[17] X. Li and K. Wang, {\em Nonparametric hypersurfaces moving by powers of Gauss Curvature,} Michigan Math. J. 66(4): 675-682 (November 2017). DOI: 10.1307/mmj/1508810813

[18] G.M. Lieberman, {\em Second Order Parabolic Differential Equations,} World Scientific, 
Hackensack, New Jersey, 2005, Reprint of the first edition (1996). 

[19] J.A. McCoy, F. Mofarreh, and G. Williams, {\em Fully nonlinear curvature flows of axially symmetric hypersurfaces with boundary 
conditions,} Annali di Matematica Pura ed Applicata, 193(5), 2013. 

[20] V. Oliker, {\em Evolution of nonparametric surfaces with speed depending on curvature, I. The 
Gauss curvature case,} Indiana Univ. Math. J. 40 (1991), no. 1, 237-258. 

[21] J. Serrin, {\em Gradient Estimates for Solutions of nonlinear elliptic and parabolic equations,} 
in Contributions to Nonlinear Functional Analysis, Edited by E.H. Zarantonello, Proceedings of a 
Symposium Conducted by the Mathematics Research Center The University of Wisconsin, April 12-14, 1971. Academic Press, New York, London, 1971. 

[22] L. Solanilla, {\em Swimming in Curved Surfaces and Gauss Curvature,} Univ. Sci. vol. 23 no.2 Bogot\'{a} May/Aug. 2018. 

[23] L. Solanilla, {\em Sobre la formulaci\'{o}n del problema de prescribir la curvatura 
de una variedad riemanniana bedimensional,} Eureka, Revista de la Licenciatura en Matem\'{a}ticas 
Aplicadas, Dic. 1998, no. 13.

[24] K. Tso, {\em Deforming a hypersurface by its Gauss-Kronecker curvature,} Comm. Pure Appl. 
Math. {\bf 38} (1985), no. 6, 867-882. 

[25] J. Urbas, {\em Complete noncompact self-similar solutions of Gauss curvature flows I. Positive 
powers,} Math. Ann. 311, 251-274 (1998). 

\vspace{1cm} 

\noindent Thalia D. Jeffres 

\noindent Department of Mathematics and Statistics 

\noindent Wichita State University 

\noindent Wichita, Kansas 

\noindent 67260-0033 

\noindent thalia.jeffres@wichita.edu 

\vspace{.5cm} 

\noindent Leonardo Solanilla 

\noindent Departamento de Matem\'{a}ticas y Estad\'{i}stica 

\noindent Universidad del Tolima 

\noindent Ibagu\'{e}, Colombia 

\noindent leonsolc@ut.edu.co

\end{document}